\documentclass[a4paper,12pt]{article}
\usepackage[english]{babel}
\usepackage[utf8]{inputenc}
\usepackage{geometry}
\usepackage{hyperref}
\usepackage{amssymb}
\usepackage{amsmath}
\usepackage{amsthm}
\usepackage{graphicx}
\usepackage{makecell}
\allowdisplaybreaks
\title{The Complex of Non-Chromatic Scales}
\author{Kathlén Kohn, Ernst Ulrich Deuker}

\theoremstyle{definition}

\theoremstyle{plain}

\theoremstyle{remark}

\begin{document}
\maketitle

\begin{abstract}
We consider the space of all musical scales with the ambition to systematize it. To do this, we pursue the idea to view certain scales as basic constituents and to ``mix'' all remaining scales from these.
The German version of this article appeared in \emph{Mitteilungen der DMV}, volume 25, issue 1. 
\end{abstract}

The musical idea of using basic constituents for the space of all scales has been suggested in the recently published book~\cite{buch} on improvisation (not only) in jazz music.
From the mathematical point of view, these constituents form a \emph{simplicial complex}, whose \emph{facets} coincide with the most widely used scales in western music~-- with the exception of the blues scale.
We will explain this connection in the following.
First, we have to clarify what exactly we mean by a scale.
The pitch space that is currently used in western music contains twelve different pitches: These are the seven \emph{natural notes} $C,D,E,F,G,A,B$ as well as the five notes $C\sharp/D\flat,D\sharp/E\flat,F\sharp/G\flat,G\sharp/A\flat,A\sharp/B\flat$,
which appear as raising/lowering of the naturals.
\begin{figure}[h!]
\centering
\includegraphics[width=0.6\linewidth]{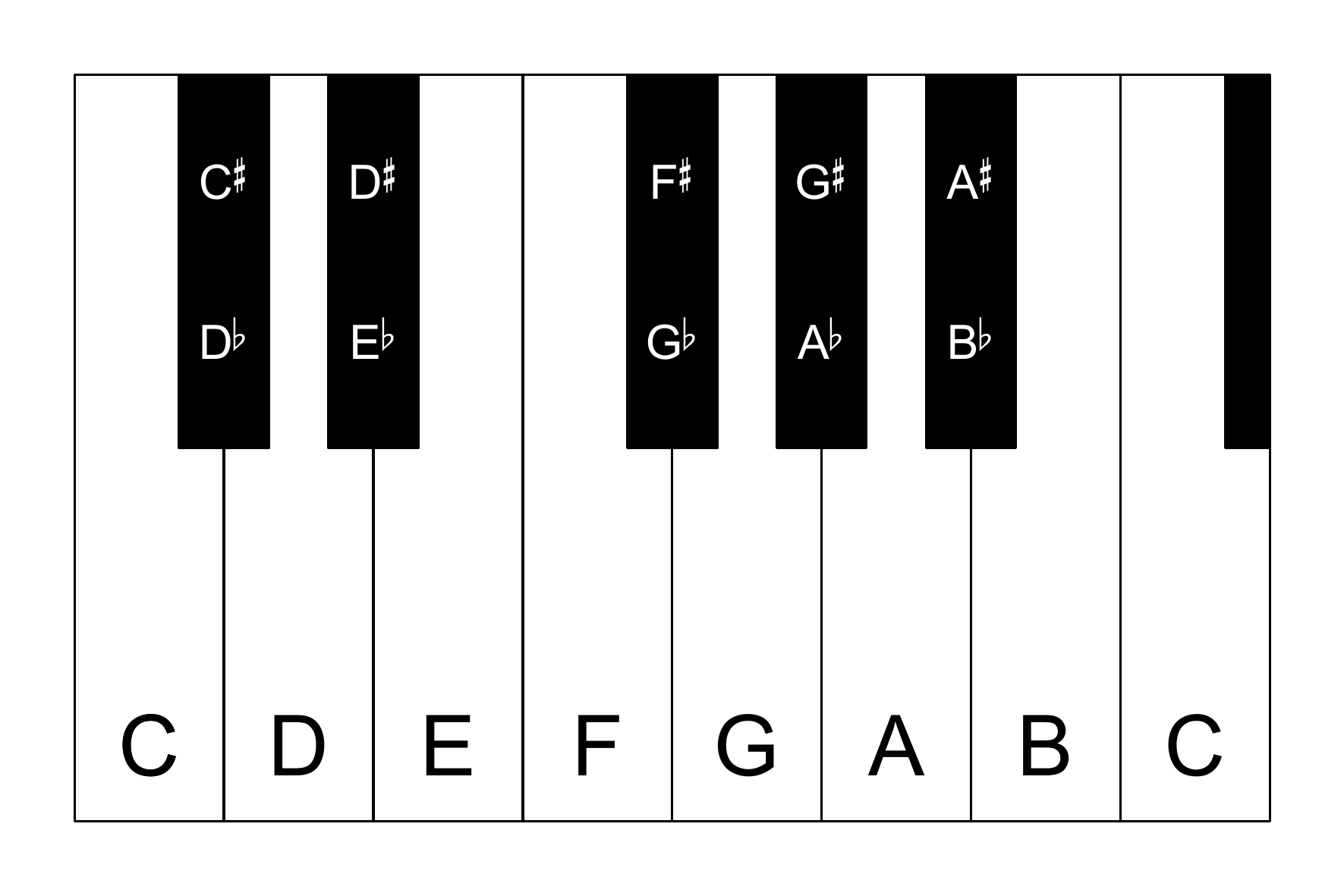}
\caption{Western pitch space in modern times with twelve pitch classes.}
\label{fig:piano}
\end{figure}

\begin{figure}[h!]
\centering
\includegraphics[width=0.6\linewidth]{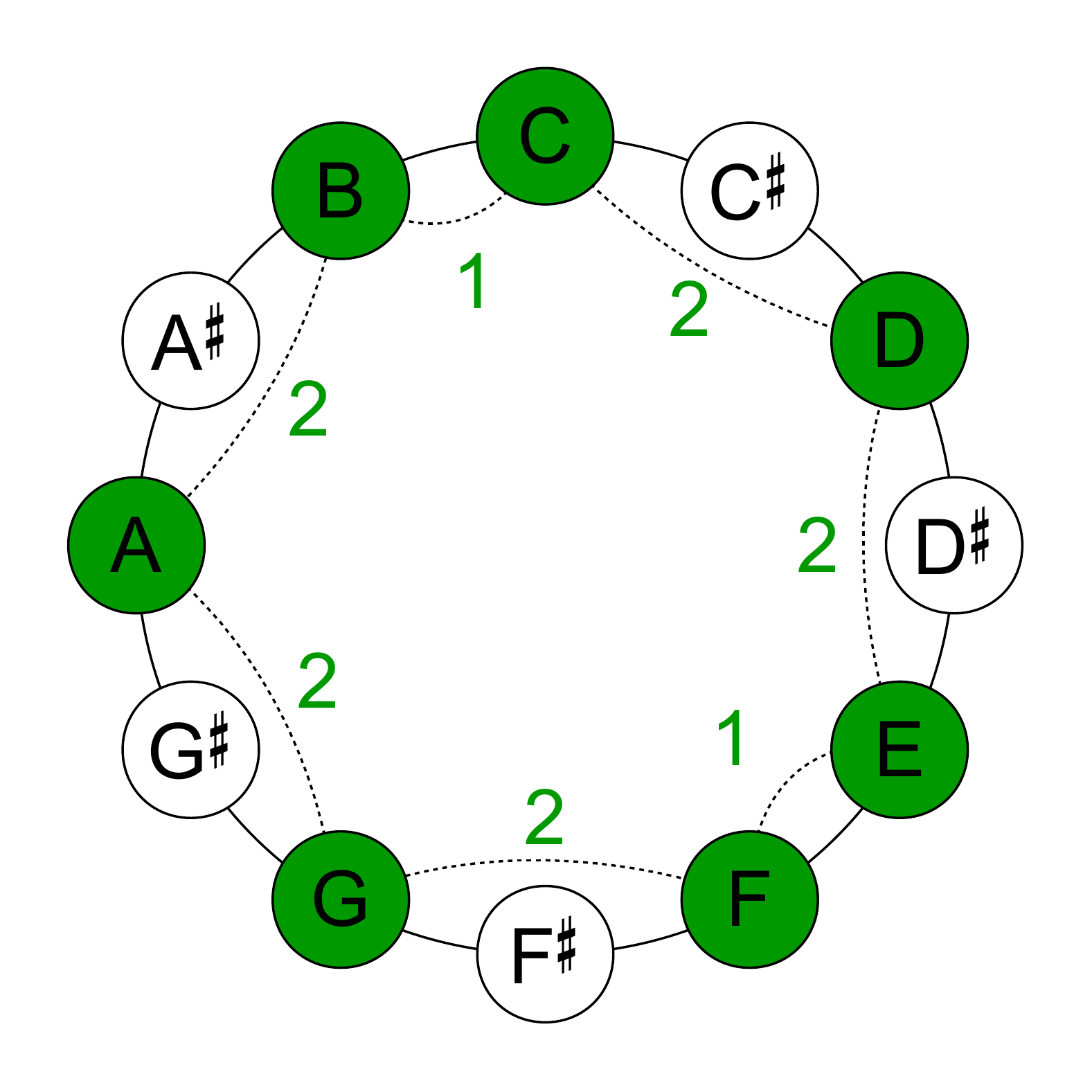}
\caption{Cyclic order of the twelve pitch classes. Marked in green: $C$ major scale with interval sequence.}
\label{fig:c}
\end{figure}

From the mathematical perspective, we do not want to distinguish the notes  $C\sharp$ and $D\flat$. In what follows, we will always use the $\sharp$-convention.
A \emph{scale} is simply defined to be a subset of  $\lbrace C, C\sharp, D, D\sharp, E, F, F\sharp, G, G\sharp, A, A\sharp, B \rbrace$.
As an example, we consider the $C$ major scale in our subset notation:
\begin{align*}
  \left\lbrace C, D, E, F, G, A, B \right\rbrace.
\end{align*}
The twelve notes have a \emph{cyclic order}. For example in the $C$ major scale we have again $C$ after the $B$ (see Fig.~\ref{fig:c}). Note that the two $C$'s in Figure~\ref{fig:piano} are different pitches, but they are said to be in the same \emph{pitch class}.
This cyclic order allows us to define a distance between two pitch classes $t_1,t_2 \in \lbrace C, C\sharp, D, D\sharp, E, F, F\sharp, G, G\sharp, A, A\sharp, B \rbrace$:
On the one hand, we can consider the clockwise distance from $t_1$ to $t_2$, and on the other hand the clockwise distance from $t_2$ to $t_1$.
The \emph{distance} between $t_1$ and $t_2$ is defined to be the minimum of both clockwise distances.
As an example we look at the pitch classes $A$ and $C$: The clockwise distance from $A$ to $C$ is three and the clockwise distance from $C$ to $A$ is nine, which means that these two pitch classes have distance three.
The sequence of distances between consecutive pitch classes in a scale in the cyclic order is called the \emph{interval sequence} of the scale.
As we can see in Figure~\ref{fig:c}, the $C$ major scale has the interval sequence 2-2-1-2-2-2-1.

The cyclic order of the twelve pitch classes implies that the interval sequence of a scale has a cyclic order, too.
For example, we can also say that the $C$ major scale has the interval sequence 2-1-2-2-2-1-2.
From the musical perspective we would say that we start the major scale from a different \emph{root} than $C$~-- in this case $D$.
This yields the widely used \emph{Gregorian modes}, \emph{$D$-Dorian} in our case.
Mathematically speaking, we would identify the scales \emph{$C$-Ionian}, \emph{$D$-Dorian}, \emph{$E$-Phrygian}, \emph{$F$-Lydian}, \emph{$G$-Mixolydian}, \emph{$A$-Aeolian} (\emph{natural $A$ minor}) and \emph{$B$-Locrian} with each other, since they all have the same interval sequence.
These modes are important in \emph{modal jazz}.
Since there are twelve major scales and thus $12 \cdot 7 = 84$ different Gregorian modes, the mentioned book by Deuker~\cite{buch} takes the didactical approach to restrict oneself to the consideration of the major scales instead of studying all Gregorian modes, although this is still common in the widespread \emph{scale theory}. A modern representative of the latter approach is for example Frank Sikora~\cite{scaleTheory}.

An interval of length one is called \emph{semitone} (or \emph{half tone}) and an interval of length two is called \emph{whole tone}.
The scale whose intervals are twelve semitones (und thus contains all twelve pitch classes) is called \emph{chromatic scale}.
Therefore, we say that a scale whose interval sequence does not contain two consecutive semitones is \emph{non-chromatic}.
This is equivalent to the scale not containing three consecutive pitch classes in the cyclic order in Figure~\ref{fig:c}.
The $C$ major scale is an example for such a non-chromatic scale. A counterexample is given by the scale
\begin{align*}
  \left\lbrace C, C\sharp, E, F, G, A, B \right\rbrace, 
\end{align*}
with interval sequence 1-3-1-2-2-2-1.
With some justification, we can view the non-chromatic scales as the ``primary colors'' in music, whereas the (partly) chromatic scales with at least two consecutive semitones could be seen as ``secondary colors''.
With this approach, the mentioned book~\cite{buch} tries to order the space of all scales systematically.

Here we consider the non-chromatic scales from the mathematical perspective.
The non-chromatic scales form a \emph{simplicial complex}.
Such a simplicial complex is defined on a ground set $\mathcal{G}$ as a set $\mathcal{K}$ of finite subsets of $\mathcal{G}$ such that, for every set $M$ in $\mathcal{K}$ and every subset $T$ of $M$, we have that $T$ is also in $\mathcal{K}$.
As an example we choose the ground set $\mathcal{G} := \lbrace 0, 1, 2 \rbrace.$
The set $\mathcal{K}_1 := \left\lbrace \lbrace 0,1 \rbrace, \lbrace 2 \rbrace \right\rbrace$ is \emph{not} a simplicial complex since $\lbrace 0 \rbrace$ is a subset of $\lbrace 0,1 \rbrace$ but not contained in $\mathcal{K}_1$.
We can extend $\mathcal{K}_1$ to a simplicial complex in the following way: 
$\mathcal{K}_2 := \left\lbrace \lbrace 0,1 \rbrace, \lbrace 0 \rbrace, \lbrace 1 \rbrace, \lbrace 2 \rbrace, \emptyset \right\rbrace$.

The set of all non-chromatic scales is a simplicial complex on the ground set of the twelve pitch classes $\mathcal{G} := \lbrace C, C\sharp, D, D\sharp, E, F, F\sharp, G, G\sharp, A, A\sharp, B \rbrace$.
Indeed, if we remove a pitch class from a non-chromatic scale, the scale stays non-chromatic.
We denote this \emph{complex of non-chromatic scales} in the following by $\mathcal{K}_{NC}$.
In other words, we can define $\mathcal{K}_{NC}$ as the set of all subsets of the twelve pitch classes in Figure~\ref{fig:c} which do not contain three consecutive pitch classes.
In the following, we want to answer three questions concerning $\mathcal{K}_{NC}$ which mathematicians typically ask about a given simplicial complex.
Moreover, we will see that these questions are musically relevant.

\section{$f$-Vector}
\begin{figure}[h!]
\centering
\includegraphics[width=0.6\linewidth]{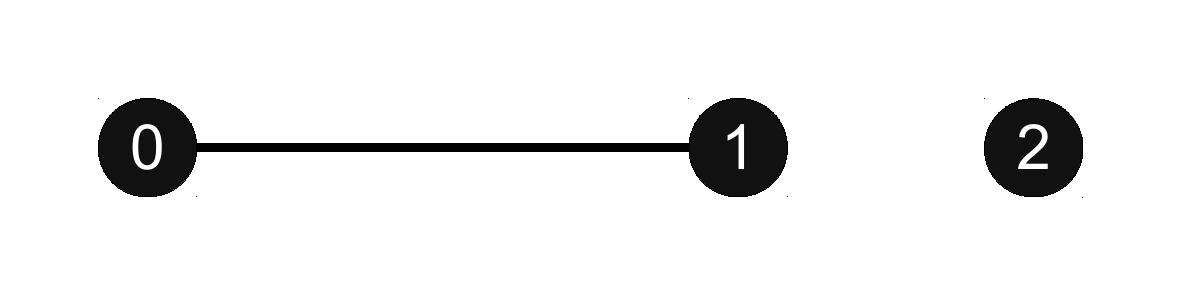}
\caption{$\mathcal{K}_2$.}
\label{fig:K2}
\end{figure}

We can represent the sets with $n$ elements in a simplicial complex geometrically as $(n-1)$-dimensional objects.
A set with one element corresponds to a point, a set with two elements to a line segment, a set with three elements to a triangle, a set with four elements to a tetrahedron, etc. 
This allows us to draw simplicial complexes (see Figure~\ref{fig:K2} for an illustration of the example complex $\mathcal{K}_2$).
Such an illustration is purely schematic; for example, embeddability of the complex is not essential.

The higher-dimensional generalization of a tetrahedron is called \emph{simplex}. 
This is where name simplicial complex stems from.
Now we consider a fixed simplicial complex and a positive integer $n$. 
We denote by $f_n$ the number of $n$-dimensional simplexes in $\mathcal{K}$ (i.e., the number of sets in $\mathcal{K}$ that have $n+1$ elements).
We view the empty set $\emptyset$, which is contained in every simplicial complex, as a $-1$-dimensional simplex and therefore we define $f_{-1} = 1$.
The list $(f_{-1},f_0,f_1, f_2, f_3, \ldots)$ of all these numbers is called \emph{$f$-vector} of $\mathcal{K}$.
The $f$-vector of the example complex $\mathcal{K}_2$ is $(1,3,1)$.

When we ask for the $f$-vector of our simplicial complex $\mathcal{K}_{NC}$, we ask how many non-chromatic scales with $0, 1, 2, 3, \ldots, 12$ pitch classes exist.
This is relevant musically; for example, one might want to know how to improvise non-chromatically.
In this case, it is particularly important to know how many pitch classes one is allowed to use maximally.
Furthermore, there are several improvisational approaches which are built up on pentatonic scales~-- these are scales with five pitch classes.
Therefore the number of non-chromatic pentatonic scales is relevant.
We calculated the $f$-vector of $\mathcal{K}_{NC}$ using the software \texttt{polymake}~\cite{DMV:polymake} (it also possible to do it by hand without too much effort):
\begin{equation*}
  \begin{array}{ccccccccc}
  (1,&12,&66,&208,&399,&456,&282,&72,&3)\\
  (f_{-1},&f_0,&f_1,&f_2,&f_3,&f_4,&f_5,&f_6,&f_7)
  \end{array}.
\end{equation*}
In particular, this implies there there is no non-chromatic scale with nine or more pitch classes.
Moreover, there are exactly three non-chromatic scales with eight pitch classes.
We will have a closer look at these three scales in the next section.
Furthermore, we note that there are 456 non-chromatic pentatonic scales.

\section{Facets}
A simplex in a simplicial complex $\mathcal{K}$ that is not contained in any other simplex of $\mathcal{K}$ is called \emph{facet}.
The example complex $\mathcal{K}_2$ has two facets, namely $\lbrace 0,1 \rbrace$ and $\lbrace 2 \rbrace$.

In musical terms, we want to know how many non-chromatic scales exist to which we cannot add a further pitch class without creating two consecutive semitones.
One example for such a scale is the $C$ major scale.
This concept of maximality is also important to improvising: It is enough to remember all maximal non-chromatic scales since all other non-chromatic scales are subsets of the maximal ones.
This is the reason why the mention book by Deuker~\cite{buch} treats these scales in detail.
He even proves in the second chapter that there are exactly 57 maximal non-chromatic scales and outlines the significance of these scales in music history.
It is remarkable that these scales are important in music as well as visible in the mathematical formulation (as facets of $\mathcal{K}_{NC}$).

In the following we describe the 57 facets of $\mathcal{K}_{NC}$.

\begin{table}[h!]
\centering
\begin{tabular}{clcl}
  \Xhline{2\arrayrulewidth}
  number of pitch classes & intervall sequence & number of scales & name\\
  \hline
 8 & 2-1-2-1-2-1-2-1 & 3 & diminished\\
 7 & 2-2-1-2-2-2-1 & 12 & major\\
 7 & 2-1-2-2-2-2-1 & 12 & melodic minor \\
 7 & 2-1-2-2-1-3-1 & 12 & harmonic minor \\
 7 & 2-2-1-2-1-3-1 & 12 & harmonic major \\
 6 & 2-2-2-2-2-2 & 2 & whole tone \\
 6 & 1-3-1-3-1-3 & 4 & augmented\\
  \Xhline{2\arrayrulewidth}
\end{tabular}
\caption{57 maximal non-chromatic scales.}
\label{table}
\end{table}

We see that the 57 maximal non-chromatic scales have only seven different interval sequences.
Looking at the first row, there are three different scales which each have eight pitch classes and the prescribed interval sequence.
\begin{figure}[h!]
\centering
\includegraphics[width=0.6\linewidth]{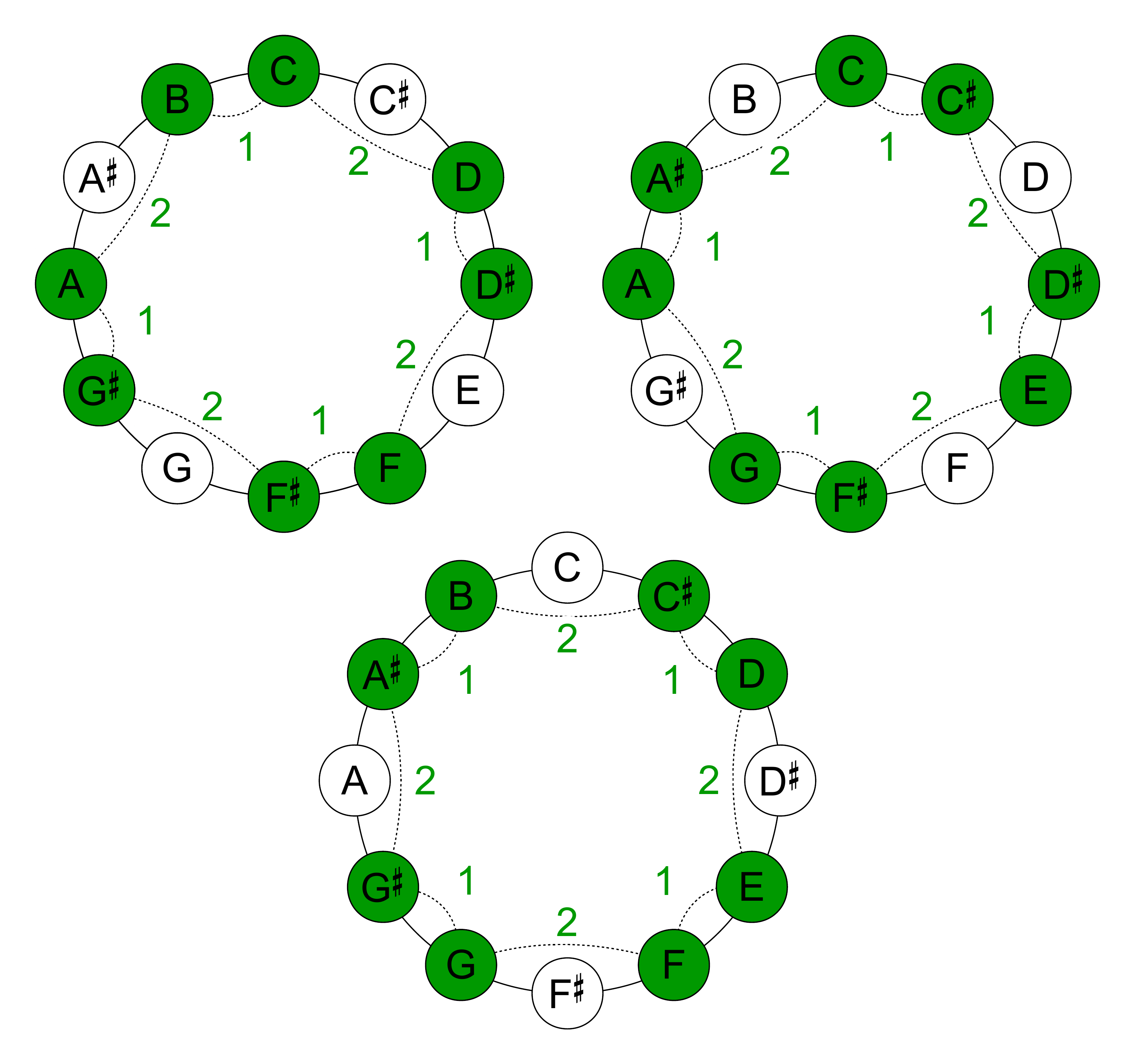}
\caption{The three diminished scales.}
\label{fig:dimin}
\end{figure}
These three scales are exactly those that we encountered when calculating the $f$-vector.
Furthermore, we see 48 maximal non-chromatic scales which consist of seven pitch classes.
These come in four different types, determined by the interval sequence.
The $f$-vector tells us that there are exactly $72=48+24$ non-chromatic scales with seven pitch classes.
The remaining $24 = 3 \cdot 8$ scales are subsets of the three scales with eight pitch classes in Figure~\ref{fig:dimin}, each of those containing eight subsets with seven elements.
Moreover, there are six maximal non-chromatic scales with six pitch classes.
In particular, we see that the facets of our simplicial complex have different dimensions.
Mathematically we say that $\mathcal{K}_{NC}$ is not \emph{pure}.

Now we want to describe the meaning of these 57 facets in music.

\textbf{2-1-2-1-2-1-2-1: }
This interval sequence has three different associated scales, see Figure~\ref{fig:dimin}.
These are called \emph{diminished scales} or even \emph{the octatonic scales}, because they are the most widely used scales among all the scales with eight pitch classes.
Nowadays they are frequently used in jazz.
There is a whole jazz-textbook devoted to these scales~\cite{diminished}.
Around the year 1900, the diminished scales were quite popular in Russia, in particular in compositions by Rimski-Korsakov.
The Dutch composer Willem Pijper used these scales extensively as well.
Thus, the diminished scales are sometimes also called Korsakovian scale or Pijper scale.

\textbf{2-2-1-2-2-2-1: }
This interval sequence corresponds to the well-known \emph{major scale} (today often called \emph{melodic major scale}). 
Our first example of a scale was the $C$ major scale, but all twelve possibilities to choose a starting note yield twelve different major scales ($D\flat$ major, $D$ major, etc.).
These are the most widely used scales in western music.
As already mentioned, the interval sequence above corresponds as well to the \emph{natural minor scales} and the other \emph{Gregorian modes}.
For example, the $C$ major scale has exactly the same pitch classes as the natural $A$ minor scale.
From the point of view of minor scales, we would write the interval sequence as 2-1-2-2-1-2-2.
This variant of minor scales it the second most common scale in occidental music, directly after the major scales.

\textbf{2-1-2-2-2-2-1: }
These are the \emph{melodic minor scales}.
There are again twelve different ones, depending on the chosen starting note.
The melodic minor scales are widespread, from early compositions~-- particularly vocal music~-- to modern pop and rock.
There is a common rule in books about the theory of harmony which says that the ascending melodic minor scale is the scale defined above, whereas the descending melodic minor scale is the same as the natural minor scale with the same staring note. 
This rule is uninteresting from the mathematical perspective, and
Deuker~\cite{buch} takes the view that
it does not make much sense from the practical musical perspective either since specific melodies are rarely linearly ascending or descending.

\textbf{2-1-2-2-1-3-1: }
This interval sequence corresponds to the \emph{harmonic minor scales}.
All twelve possible starting notes give twelve different such scales.
In contrast to melodic minor scales, the harmonic minor scales are mainly used to build chords, 
but in particular in compositions by Mozart or Schubert they do appear in melodies as well.

\textbf{2-2-1-2-1-3-1: }
This last type of seven-note maximal non-chromatic scales has again twelve different representatives.
These are called \emph{harmonic major scales}.
They are conceived as a mixture of major and harmonic minor scales, and they are found most frequently in jazz music.

\textbf{2-2-2-2-2-2: }
This is the \emph{whole tone scale}.
Due to its high symmetry, there are just two different such scales, depending on the starting note.
Since there are no semitones at all, the choice of the starting note is a relatively arbitrary determination and there is no audible root in the scale.
As a consequence, the impression of the whole tone scale is often described as levitating.
In modern western music (except jazz), this scale is not used that often, but in impressionism~-- in particular in compositions by Debussy~-- it played a major role, and it was applied before this era by Franz Liszt and Rimski-Korsakow.

\textbf{1-3-1-3-1-3: }
The \emph{augmented scales} belong to this interval sequence.
Analogously to the case of the diminished scales (see Fig.~\ref{fig:dimin}), there are four different augmented scales.
These appeared already in compositions by Franz Liszt, but they were used increasingly in the 20th century by composers like B\'ela Bart\'ok or Arnold Sch\"onberg as well as in jazz.
There is even a jazz-textbook which is explicitly dedicated to the augmented scales~\cite{augmented}.

At the end of this section, we want to deduce why there are exactly the seven possibilities in Table~\ref{table} for an interval sequence $\mathcal{I}$ of a maximal non-chromatic scale.

\textbf{Observation 1:}
$\mathcal{I}$ cannot contain an interval of length four or more.
If $\mathcal{I}$ would have such an interval~--e.g., between the pitch classes $C$ and $E$~--, then we could add a pitch class in the middle of the interval ($D$ in the example) and the scale would stay non-chromatic.

\textbf{Observation 2:}
If $\mathcal{I}$ contains an interval of length three, then the right and the left neighbor of this interval in $\mathcal{I}$ must be both semitones. In other words, if 3 is contained in $\mathcal{I}$, then $\mathcal{I}$ has to contain the sequence 1-3-1.
Indeed, if $\mathcal{I}$ would contain the sequence 3-2~-- e.g., the pitch classes $C$, $D\sharp$ and $F$ are in the scale~--, then we could add a pitch class from the interval of length three (in the example $D$) to the scale such that the scale would stay non-chromatic.

Using these two observations, we can distinguish a few cases to explore all seven possibilities for $\mathcal{I}$.

\textbf{Case 1:}
$\mathcal{I}$ contains a 3 and thus~-- by Observation~2~-- the sequence 1-3-1.
Either we extend this sequence with a 3 or with a 2.
If we choose to extend it with a 3, we see that $\mathcal{I}$ has to contain the sequence 1-3-1-3-1.
Since all numbers in $\mathcal{I}$ have to add up to twelve, there is just one possibility to extend the sequence to an interval sequence of a non-chromatic scale:
\begin{align*}
  1-3-1-3-1-3.
\end{align*}
If we would have chosen to extend the sequence 1-3-1 on the left side with a 3, the same interval sequence would have been our result.
Hence, we have just one case left to consider, namely that $\mathcal{I}$ contains the sequence 2-1-3-1-2. 
By Observation~2, we have to complete this sequence with a 2 and a 1, which yields two possibilities:
\begin{align*}
  2-2-1-2-1-3-1,\\
  2-1-2-2-1-3-1.
\end{align*}

\textbf{Case 2:}
$\mathcal{I}$ contains no 3, and as a consequence, consists only of ones and twos.
Since the sum of all numbers in $\mathcal{I}$ has to be twelve, the number of ones in $\mathcal{I}$ is even.
If $\mathcal{I}$ contains no ones at all, we know that $\mathcal{I}$ looks as follows:
\begin{align*}
  2-2-2-2-2-2.
\end{align*}
If $\mathcal{I}$ contains exactly two ones, then there are only two possibilities to arrange them:
\begin{align*}
  2-1-2-2-2-2-1,\\
  2-2-1-2-2-2-1.
\end{align*}
Due to the cyclic order, the further possibilities 2-2-2-1-2-2-1 and 2-2-2-2-1-2-1 are in fact identically equal to the two interval sequences above.
If $\mathcal{I}$ contains four ones, it is already uniquely determined:
\begin{align*}
  2-1-2-1-2-1-2-1.
\end{align*}
Moreover, this shows that $\mathcal{I}$ cannot contain more than four ones.

\section{Topology}

\emph{Topology} is a branch of mathematics that investigates which properties of spaces are preserved under continuous deformations.
We say that a deformation is \emph{continuous} if~-- roughly speaking~-- it does not require one to cut or glue parts of the space.
Typical examples of continuous deformations are stretching and bending of spaces.
For example, a ball can be continuously deformed to a cube.
Since this deformation can also be reversed continuously, we say that ball and cube are \emph{homeomorphic}.
A mug with a handle is homeomorphic to a donut with a hole.
However, such a donut is not homeomorphic to a ball, because its hole cannot be filled by a continuous deformation. 
Thus, the number of \emph{holes} is an example of a property that homeomorphic spaces have in common.

The idea of homeomorphy can be extended to \emph{homotopy equivalence}, which also preserves the number of holes.
Let $X$ and $Y$ be two spaces with two continuous maps $g: X \to Y$ and $f:Y \to X$.
If the compositions $f \circ g: X \to X$ and $g \circ f: Y \to Y$ are the identities on $X$ and $Y$, respectively (i.e., $f(g(x)) = x$ for all $x \in X$ and $g(f(y))  = y$ for all $y \in Y$), then $X$ and $Y$ are homeomorphic.
If the map $f \circ g$ can be deformed continuously to the identity on $X$ and if $g \circ f$ can be deformed continuously to the identity on $Y$, then we say that $X$ and $Y$ are \emph{homotopy equivalent}.
Hence, every homeomorphism is a homotopy equivalence.
Examples of homotopy equivalences that are not homeomorphisms are to thicken or to squash spaces.
For instance, a ball and a point are homotopy equivalent but not homeomorphic.

\subsection{Holes in Simplicial Complexes}
We can also count holes in a given simplicial complex $\mathcal{K}$.
For this, we determine the \emph{(reduced) simplicial homology} of the complex.
We give here a short introduction to this topic and refer the interested reader to~\cite{topo}.
As in the definition of the $f$-vector, we consider all $n$-dimensional simplexes in $\mathcal{K}$.
A \emph{formal sum} over these simplexes denotes simply a sum of the simplexes with rational coefficients, which we write down without evaluating it.
We denote by $C_n(\mathcal{K})$ the set of all such formal sums.
This is a vector space over the rational numbers $\mathbb{Q}$.
In our example complex $\mathcal{K}_2$, we have among others the following formal sums:
\begin{equation*}
  \frac{3}{2} \lbrace  0,1 \rbrace \in C_1(\mathcal{K}_2), \;
  2 \lbrace 0 \rbrace - \lbrace 2 \rbrace \in C_0(\mathcal{K}_2), \;
  - \emptyset \in C_{-1}(\mathcal{K}_2).
\end{equation*}

The \emph{boundary operator} maps an $n$-dimensional simplex to its $(n-1)$-dimensional boundary.
It maps, for instance, a triangle to its three bounding line segments, and a tetrahedron to its four faces (which are triangles).
Formally, the boundary operator $\partial_n: C_n (\mathcal{K}) \to C_{n-1}(\mathcal{K})$ is defined by mapping every set in $\mathcal{K}$ with $(n+1)$ elements to an alternating sum of its subsets with $n$ elements.
For this, the elements of a set with $(n+1)$ elements will be omitted one after the other, in an order that was fixed beforehand.
In our example $\mathcal{K}_2$ this means:
\begin{align*}
  \partial_1(\lbrace 0,1 \rbrace) = \lbrace 1 \rbrace - \lbrace 0 \rbrace,\\
  \partial_0(\lbrace 0 \rbrace) = \partial_0(\lbrace 1 \rbrace) = \partial_0 (\lbrace 2 \rbrace) = \emptyset.
\end{align*}
The formal sums in $C_n(\mathcal{K})$ that get mapped to 0 by the boundary operator are called \emph{$n$-cycles}.
The set of all these cycles is a vector space, which we denote by $Z_n(\mathcal{K})$.
We call the image $\partial_{n+1}(\Sigma)$ of a formal sum $\Sigma \in C_{n+1}(\mathcal{K})$ an \emph{$n$-boundary}.
The set of all these boundaries is also a vector space, denoted $R_n(\mathcal{K})$.
Now one can verify immediately that every $n$-boundary is an $n$-cycle, i.e., $R_n(\mathcal{K}) \subseteq Z_n(\mathcal{K})$.
We say that two $n$-cycle are \emph{equivalent} if their difference is an $n$-boundary.
The \emph{equivalence class} of an $n$-cycle $\zeta$ is the set of all $n$-cycles that are equivalent to  $\zeta$.
We denote by $H_n(\mathcal{K})$ the set of all equivalence classes.
This is again a vector space, called the \emph{$n$th homology group} of $\mathcal{K}$.

In our example $\mathcal{K}_2$, the 0-cycles are all formal sums of the form
$a \lbrace 0 \rbrace + b \lbrace 1 \rbrace + c \lbrace 2 \rbrace$, where $a,b,c$ are rational numbers with $a+b+c=0$.
The 0-boundaries look as follows: $d \lbrace 0 \rbrace - d \lbrace 1 \rbrace$ with $d \in \mathbb{Q}$.
Thus, the equivalence class of a 0-cycle $a \lbrace 0 \rbrace + b \lbrace 1 \rbrace + c \lbrace 2 \rbrace$ consists of all 0-cycles that have the same coefficient $c$ in front of $\lbrace 2 \rbrace$.
Hence, such an equivalence class in uniquely determined by the coefficient $c$.
Since $c$ is an arbitrary rational number, we say that $H_0(\mathcal{K}_2)$ is \emph{isomorphic} to $\mathbb{Q}$.
The only 1-cycle of $\mathcal{K}_2$ is 0 (i.e., the formal sum with coefficients 0). 
Due to the fact that $\mathcal{K}_2$ contains no set with 3 elements, we can say that 0 is also the only 1-boundary of $\mathcal{K}_2$.
This shows that $H_1(\mathcal{K}_2)$ contains nothing but 0.

When considering the illustration of $\mathcal{K}_2$ in Figure~\ref{fig:K2}, we see that $\mathcal{K}_2$ is not \emph{connected}.
This means that one cannot draw $\mathcal{K}_2$ without lifting the pen.
Instead, $\mathcal{K}_2$ has two connected components. 
We can read off the number of these componenents from the first homology group: 
for every simplicial complex $\mathcal{K}$, the vector space $H_0(\mathcal{K})$ is isomorphic to the vector space $\mathbb{Q}^{k-1}$, where $k$ is the number of connected components of $\mathcal{K}$.
We use the convention that $\mathbb{Q}^0 = \lbrace 0 \rbrace$.

We cannot only count the connected components of a simplicial complex with homology, but also the number of holes.
To explain this, we consider a filled and an unfilled triangle as example complexes: 
\begin{align*}
  \mathcal{K}_\vartriangle &:= \lbrace \lbrace 0,1 \rbrace, \lbrace 1,2 \rbrace, \lbrace 0,2 \rbrace, \lbrace 0 \rbrace, \lbrace 1 \rbrace, \lbrace 2 \rbrace, \emptyset  \rbrace,\\
  \mathcal{K}_\blacktriangle &:= \lbrace \lbrace 0,1,2 \rbrace, \lbrace 0,1 \rbrace, \lbrace 1,2 \rbrace, \lbrace 0,2 \rbrace, \lbrace 0 \rbrace, \lbrace 1 \rbrace, \lbrace 2 \rbrace, \emptyset  \rbrace.
\end{align*}
The illustration of both complexes looks as follows:
\begin{figure}[h!]
\centering
\includegraphics[width=0.6\linewidth]{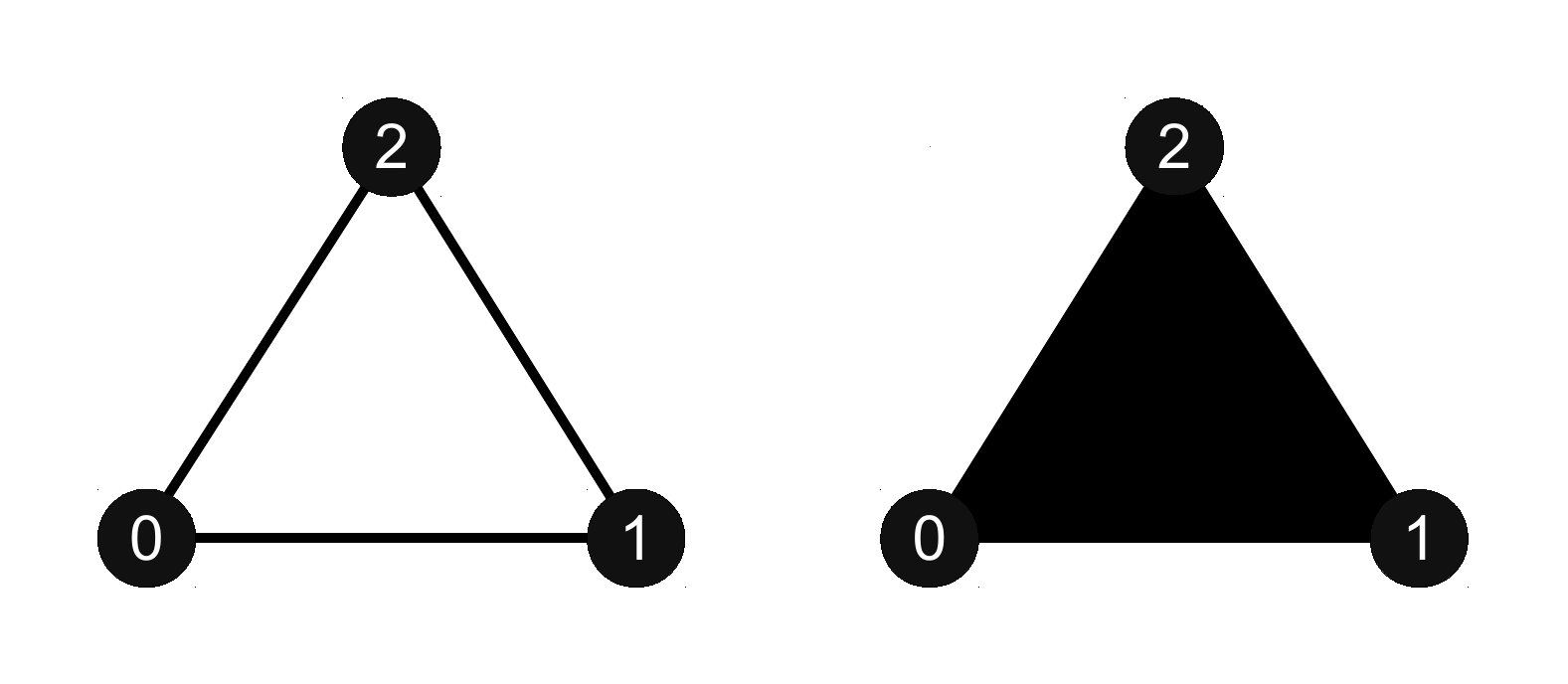}
\caption{$\mathcal{K}_\vartriangle$ and $\mathcal{K}_\blacktriangle$.}
\label{fig:triangles}
\end{figure}

Since both triangles are connected, we know that $H_0(\mathcal{K}_\vartriangle)$ and $H_0(\mathcal{K}_\blacktriangle)$ must be isomorphic to $\lbrace 0 \rbrace$.
We can verify this directly:
Both for $\mathcal{K}_\vartriangle$ and for $\mathcal{K}_\blacktriangle$, the first two boundary operators are equal:
\begin{align*}
  \partial_1(\lbrace 0,1 \rbrace) = \lbrace 1 \rbrace - \lbrace 0 \rbrace,\\
  \partial_1(\lbrace 1,2 \rbrace) = \lbrace 2 \rbrace - \lbrace 1 \rbrace,\\
  \partial_1(\lbrace 0,2 \rbrace) = \lbrace 2 \rbrace - \lbrace 0 \rbrace,\\
  \partial_0(\lbrace 0 \rbrace) = \partial_0(\lbrace 1 \rbrace) = \partial_0 (\lbrace 2 \rbrace) = \emptyset.
\end{align*}
Thus, all 0-cycles of both simplicial complexes are of the form
$a \lbrace 0 \rbrace + b \lbrace 1 \rbrace + c \lbrace 2 \rbrace$, where $a,b,c$ are rational numbers with $a+b+c=0$.
The 0-boundaries look like
\begin{align*}
d(\lbrace 1 \rbrace - \lbrace 0 \rbrace)+e(\lbrace 2 \rbrace - \lbrace 1 \rbrace)+f(\lbrace 2 \rbrace - \lbrace 0 \rbrace) \\ = 
(-d-f) \lbrace 0 \rbrace + (d-e) \lbrace 1 \rbrace + (e+f) \lbrace 2 \rbrace
\end{align*}
with $d,e,f \in \mathbb{Q}$.
From this we deduce that every 0-cycle is a 0-boundary.
Hence, there is just one equivalence class of 0-cycles, and $H_0(\mathcal{K}_\vartriangle) = H_0(\mathcal{K}_\blacktriangle)$ is a set with one element, which is consequently isomorphic to $\lbrace 0 \rbrace$, as claimed above.

The 1-cycles of $\mathcal{K}_\vartriangle$ and $\mathcal{K}_\blacktriangle$ are of the form 
$d \lbrace 0,1 \rbrace + d \lbrace 1,2 \rbrace - d \lbrace 0,2 \rbrace$.
The boundary operator $\partial_2$ of the filled triangle $\mathcal{K}_\blacktriangle$ maps $\lbrace 0,1,2 \rbrace$ to
$\lbrace 1,2 \rbrace - \lbrace 0,2 \rbrace + \lbrace 0,1 \rbrace$.
That is why the 1-boundaries of $\mathcal{K}_\blacktriangle$ look exactly like the 1-cycles, and $H_1(\mathcal{K}_\blacktriangle)$ contains just one equivalence class.
We conclude again that $H_1(\mathcal{K}_\blacktriangle)$ is isomorphic to $\mathbb{Q}^0$.
Since the unfilled triangle $\mathcal{K}_\vartriangle$ contains no sets with three elements, the only 1-boundary of $\mathcal{K}_\vartriangle$ is 0.
Therefore, different 1-cycles of $\mathcal{K}_\vartriangle$ can never be equivalent.
We get that $H_1(\mathcal{K}_\vartriangle)$ is isomorphic to $\mathbb{Q}^1$, because the coefficient $d$ of the 1-cycles is an arbitrary rational number.
The triangle $\mathcal{K}_\vartriangle$ has a (one-dimensional)hole, whereas $\mathcal{K}_\blacktriangle$ has no holes.
Thus, we have for both triangle that the dimension of $H_1(\mathcal{K}_\vartriangle)$ or $H_1(\mathcal{K}_\blacktriangle)$, respectively, is equal to the number of one-dimensional holes of $\mathcal{K}_\vartriangle$ or $\mathcal{K}_\blacktriangle$, respectively.

Furthermore, we see in Figure~\ref{fig:K2} that $\mathcal{K}_2$ has no one-dimensional holes, 
which fits to our calculation that $H_1(\mathcal{K}_2) = \lbrace 0 \rbrace$.
We can generalize this idea:
Intuitively, the dimension of $H_n(\mathcal{K})$ counts the $n$-dimensional holes in a simplicial complex $\mathcal{K}$.

\subsection{Holes in $\mathcal{K}_{NC}$}
Now one might ask the question how many holes of which dimensions are in the complex of non-chromatic scales, and what is the musical meaning of these holes?
Using the software \texttt{polymake}, we obtain that $H_5(\mathcal{K}_{NC})$ is isomorphic to $\mathbb{Q}^3$ and that all other homology groups are simply $\lbrace 0 \rbrace$.
As a consequence we can conclude that $\mathcal{K}_{NC}$ has three 5-dimensional holes.
Moreover, the software \texttt{polymake} gives us a basis for the homology for each of the three holes, i.e., a set of hexatonic (6-note) scales that defines the boundary of the hole.
Furthermore, we can check with \texttt{polymake} that the union of those hexatonic scales that define the boundary of hole is homeomorphic to the \emph{5-sphere}.
The \emph{$n$-sphere} is the boundary of an $(n+1)$-dimensional ball.
For instance, the unfilled triangle in Figure~\ref{fig:triangles} is homeomorphic to the 1-sphere, which is simply the edge of a circle.
Roughly speaking, we can imagine the complex of non-chromatic space to be composed of three 5-spheres.

The 5-sphere consists~-- as mentioned above~-- of hexatonic scales.
Additionally, all subscales of such a hexatonic scale lie also on the respective sphere; just as the unfilled triangle in Figure~\ref{fig:triangles} consists of sets with two elements, but their subsets with one element lie as points on the triangle.
Now we still have scales with seven and eight pitch classes in $\mathcal{K}_{NC}$.
We explain why these do not play any role in the topological picture of the three spheres with so-called \emph{collapses} (see e.g.~\cite{kozlov}).
Whenever a simplicial complex $\mathcal{K}$ has a facet $F$ with a subset $M \subseteq F$ such that the following two properties hold:
\begin{enumerate}
  \item $M$ contains exactly one element less than $F$,
  \item $F$ is the only facet that contains $M$,
\end{enumerate}
then we can remove $F$ and $M$ from the complex $\mathcal{K}$ and we obtain another simplicial complex.
This process is referred to as a collapse and it is an example for a homotopy equivalence of simplicial complexes.
In particular, the homology groups of the smaller complex are isomorphic to the homology groups of $\mathcal{K}$, and both simplicial complexes have the same number of holes.
Let us observe the following example to understand better how a collapse looks and why it leaves the number of holes unchanged:
\begin{figure}[h!]
\centering
\includegraphics[width=0.6\linewidth]{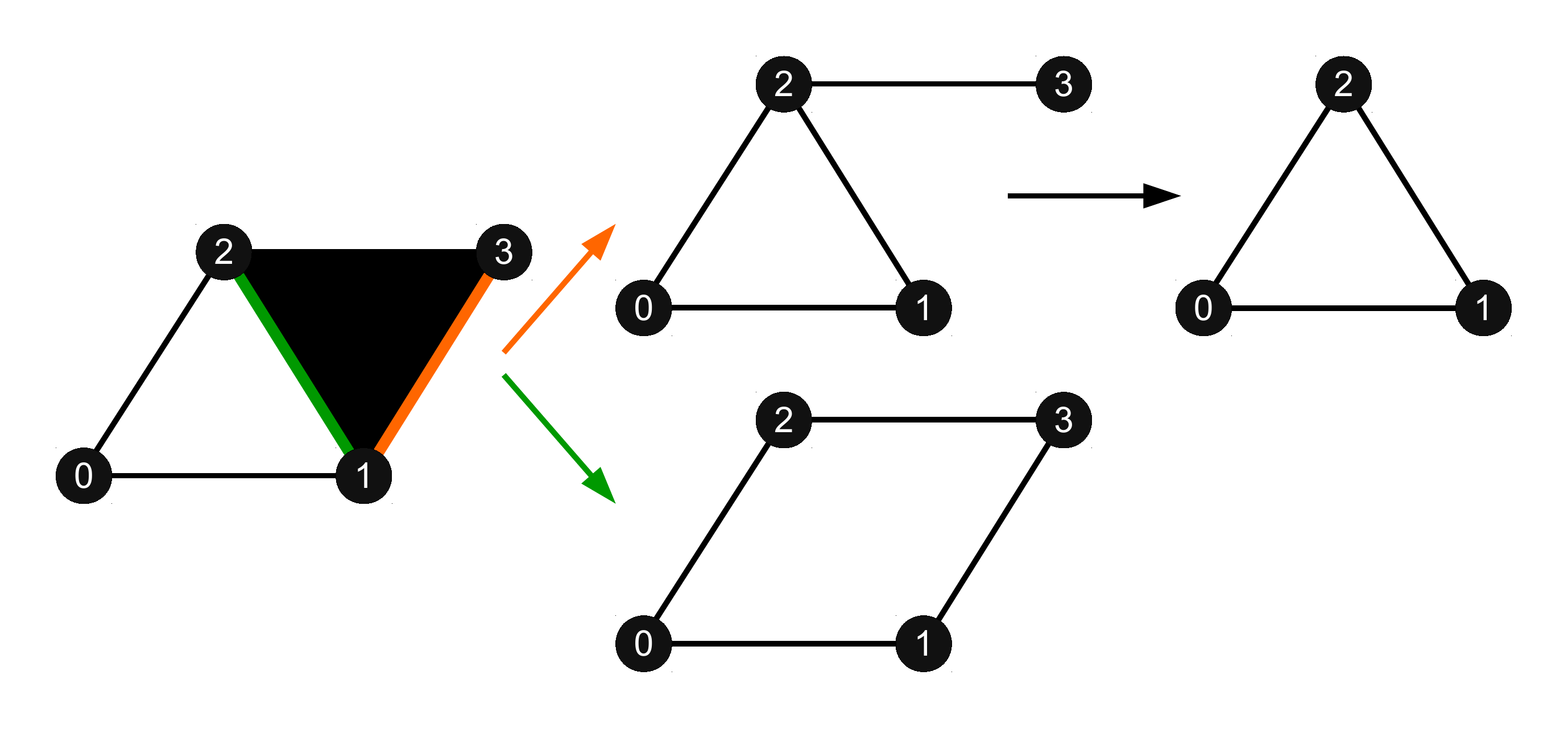}
\caption{The left simplicial complex has 3 facets, but only one of these can be collapsed: $\lbrace 1,2,3 \rbrace$. 
In doing so, we can either remove the inner edge (green) or one of the two outer edges (e.g., orange). 
In the first case, we cannot perform a second collapse after the first one.
In the second case, there is exactly one further collapse possible.}
\label{fig:collapse}
\end{figure}

In our complex of non-chromatic scales, we remove successively all scales with seven and eight pitch classes by using collapses.
For example, if we remove the pitch class $C$ from the first eight-note scale in Figure~\ref{fig:dimin}, then we obtain a seven-note scale which can be extended to an eight-note non-chromatic scale in exactly one way.
Thus, we can perform a collapse with these two scales.
We consider the $C$ major scale as a second example.
Here we can remove the pitch class $D$ such that the resulting hexatonic scale is contained in only one maximal non-chromatic scale.
Therefore, we can collapse the $C$ major scale and its hexatonic subscale.
In this manner, we can remove all scales with seven and eight pitch classes from $\mathcal{K}_{NC}$.
This means that there are higher-dimensional scales (namely with seven and eight pitch classes) glued onto the three 5-spheres in $\mathcal{K}_{NC}$ (which have hexatonic scales in their boundaries), just as the filled triangle is glued onto the unfilled triangle in Figure~\ref{fig:collapse}.
By the way, it is impossible to collapse the hexatonic facets of $\mathcal{K}_{NC}$.

We are left with the question how our three spheres look exactly.
In particular, we want to understand which hexatonic scales define the boundary of which sphere.
For this, we consider the \emph{augmented triad} $C-E-G\sharp$.
After removing this triad from our twelve pitch classes, we obtain the following nine-note scale which was used by the French composer Messiaen:
\begin{figure}[h!]
\centering
\includegraphics[width=0.6\linewidth]{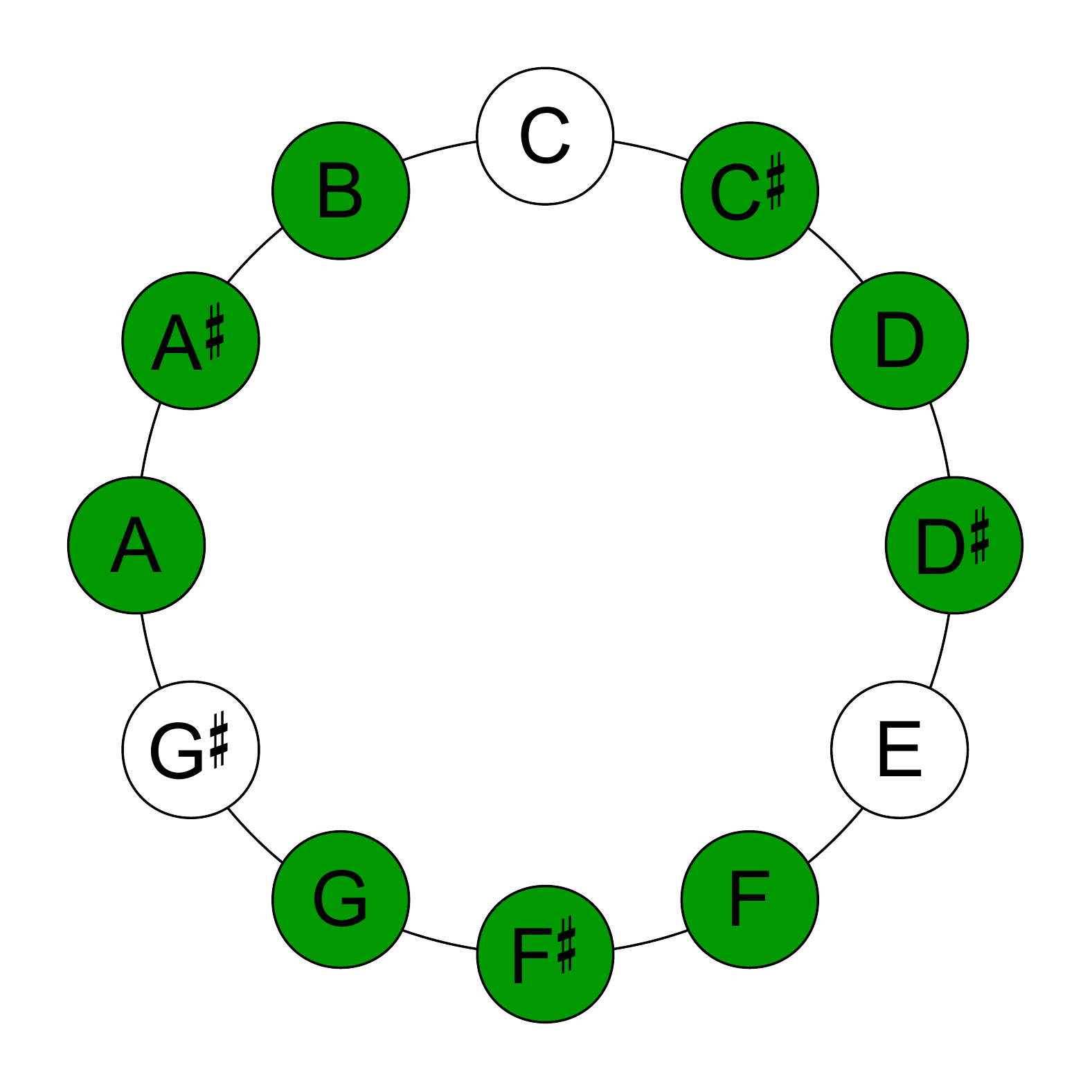}
\caption{Messiaen's nine-note scale.}
\label{fig:mess}
\end{figure}

This scale does not contain non-chromatic scales with seven or eight pitch classes, but 27 non-chromatic scales with six pitch classes.
The latter scales can be obtained as follows:
choose two pitch classes out of $C\sharp - D - D\sharp$ (there are three possibilities to do this), then pick two pitch classes from $F - F\sharp - G$ (again, there are three possibilities), and finally choose two pitch classes out of $A - A\sharp - B$ (which gives three more possibilities).
Our analysis with \texttt{polymake} shows that these 27 hexatonic scales (and their subscales) form one of our three 5-spheres.

\begin{figure}[h!]
\centering
\includegraphics[width=0.6\linewidth]{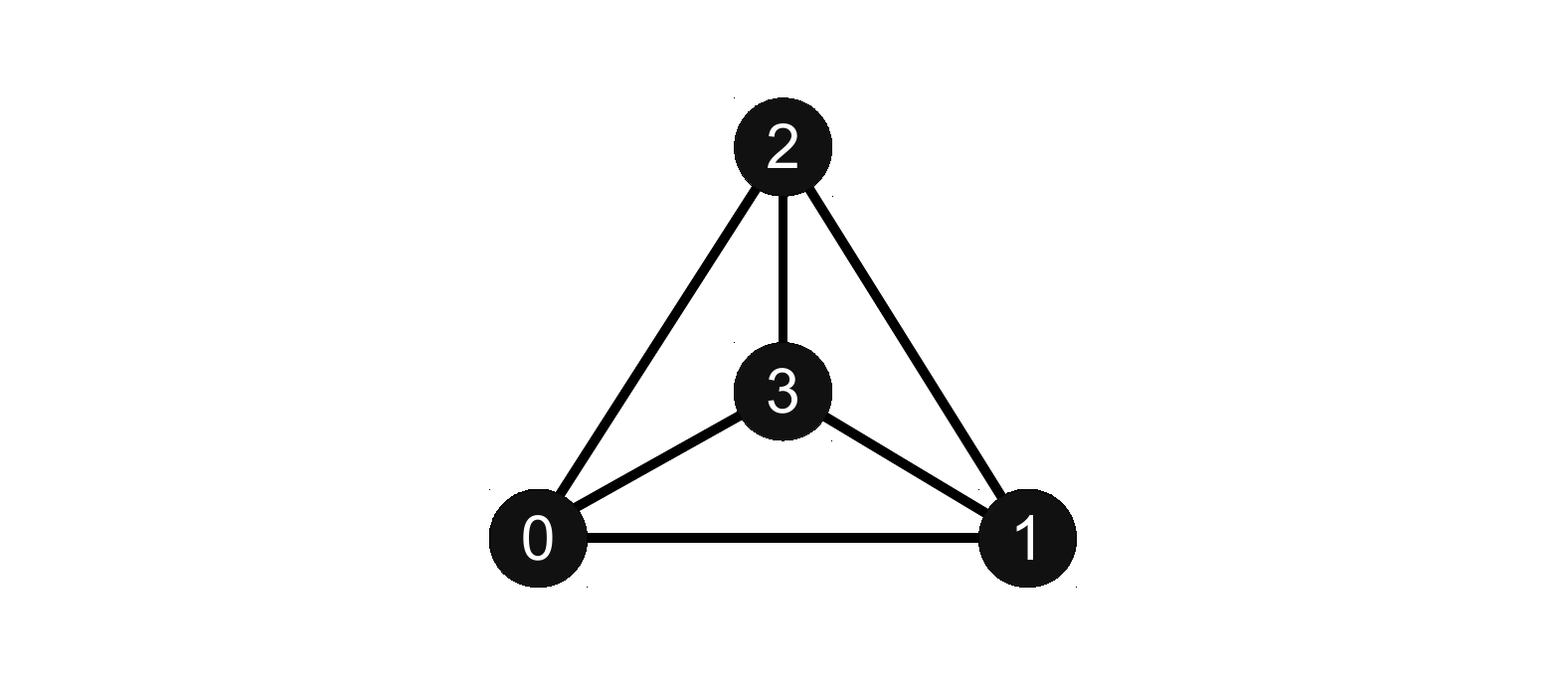}
\caption{simplicial complex with facets $\lbrace 0,1 \rbrace$, $\lbrace 0,2 \rbrace$, $\lbrace 0,3 \rbrace$, $\lbrace 1,2 \rbrace$, $\lbrace 1,3 \rbrace$, $\lbrace 2,3 \rbrace$.}
\label{fig:tetra}
\end{figure}

As in the case of the three different diminished scales in Figure~\ref{fig:dimin}, we see that there are four different Messiaen scales with nine pitch classes, which are obtained by omitting an augmented triad.
Each of these four scales contains 27 non-chromatic hexatonic scales, which form a sphere in $\mathcal{K}_{NC}$.
However, we claimed before that there are only three holes in $\mathcal{K}_{NC}$. 
The reason for this is that the four Messiaen spheres described above are not independent of each other in $\mathcal{K}_{NC}$.
This principle can be easily explained by looking at the skeleton of a tetrahedron:
Whenever we draw the simplicial complex in Figure~\ref{fig:tetra}, we see only three one-dimensional holes, although there are four 1-spheres (unfilled triangles).

\begin{figure}[h!]
\centering
\includegraphics[width=0.6\linewidth]{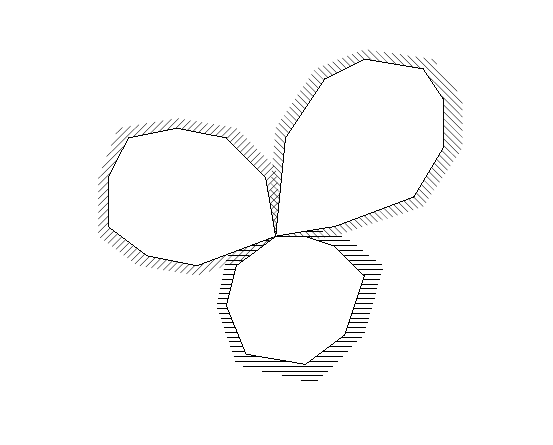}
\caption{Schematic illustration of the topology of $\mathcal{K}_{NC}$.}
\label{fig:balloons}
\end{figure}

\begin{figure}[h!]
\centering
\includegraphics[width=\linewidth]{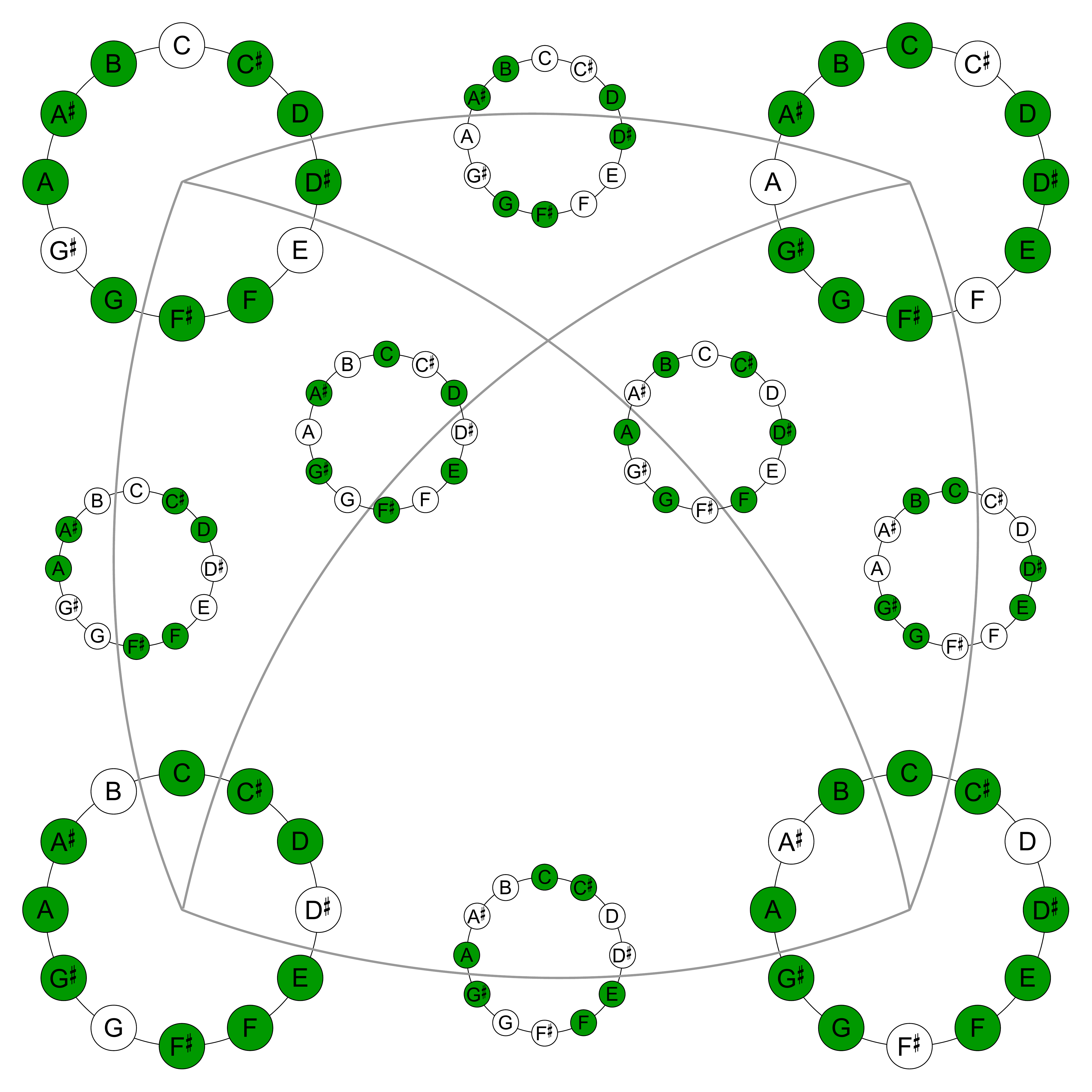}
\caption{The hexatonic facets of $\mathcal{K}_{NC}$ are the pairwise intersections of the four Messiaen scales.}
\label{fig:4messiaen}
\end{figure}

One can show that any three of the four Messiaen spheres are indeed independent of each other in $\mathcal{K}_{NC}$, i.e., that they form a basis of the vector space $H_5(\mathcal{K}_{NC})$.
Let us finally observe that each 5-sphere described above contains exactly three of the hexatonic facets of $\mathcal{K}_{NC}$.
In fact, any two distinct Messiaen scales intersect in a hexatonic facet of $\mathcal{K}_{NC}$, and all six hexatonic facets come as such an intersection (see Figure~\ref{fig:4messiaen}).
Three distinct Messiaen scales intersect in an augmented triad, i.e., in a filled triangle as in Figure~\ref{fig:triangles}.
In particular, we can successively perform collapses on the intersection of any two Messiaen spheres to obtain the common augmented triad, and afterwards we can collapse this triad until there is just a point left.
That is why we can say from a topological perspective that the three spheres are glued together at a point.
Hence, we can think of $\mathcal{K}_{NC}$ as three 6-dimensional ballons, which are held together at one point and have a higher-dimensional layer glued onto them (see Figure~\ref{fig:balloons}). 

We have described the 57 facets of the simplicial complex $\mathcal{K}_{NC}$ as well as its topological structure from the point of view of mathematics and music.
In \emph{Mitteilungen der DMV}, volume 22, issue 4, there is an article on ``Mathematik \emph{und} Musik?\grqq'' written by the mathematician and pianist Christian Krattenthaler, which takes the view that mathematics and music are actually irrelevant to each other.
However, we came to the conclusion that this is not always the case, since the simplicial complex $\mathcal{K}_{NC}$ has both mathematical and musical meaning, and there are many more connections between musical practice (e.g., change of scales) and mathematics that can be investigated. 

We thank Michael Joswig for helpful discussions, suggestions and ideas,
Bernd Sturmfels, who introduced the two authors to each other, 
and Madeline Brandt for revising the English version of this article.

\bibliographystyle{alpha}
\bibliography{literaturENG}

\end{document}